\documentclass[a4paper]{article}
\usepackage{geometry,amssymb,hyperref}
\usepackage[many]{tcolorbox}

\newtheorem{te}{Theorem}

\newtheorem{definicja}{Definition}
\newenvironment{de}{\begin{definicja}}{$\blacktriangle$ \end{definicja}}

\newenvironment{dkz}{Proof:}{$\square$}
\newcommand{\torka}[1]{\langle{#1}\rangle}

\newcommand{\upar}[2]{\torka{{#1},{#2}}}
\newcommand{\izbaci}[1]{}

\newcommand{\Part}[1]{{\cal P}{(#1)}}
\newcommand{\dom}[1]{\mathrm{dom}({#1})}

\newcommand{\funure}[3]{{\mathrm{Fn}}({#1},{#2},{#3})}

\newcommand{\ZFC}{\mathrm{ZFC}}

\newcommand{\skupsvih}[2]{{\{{#1}\,|\,{#2}\}}}
\newcommand{\jedef}{:=}
\newcommand{\ekdef}{\stackrel{\mathrm{def}}\Leftrightarrow}

\newcommand{\cf}{{\mathrm{cf}}}
\newcommand{\HLT}{Halpern and L\" auchli's theorem}
\author{Nedeljko Stefanovi\' c}
\title{The coloring principle for the product of Polish spaces and the \HLT{}}
{\begin{tcolorbox}[breakable, title=Removed, colback=yellow]}%
	{\end{tcolorbox}}
{\begin{tcolorbox}[breakable, enhanced, title=Added, colback=green, before upper={\parindent15pt}]}%
	{\end{tcolorbox}}
{\begin{tcolorbox}[breakable, title=Changed, colback=red]}%
	{\end{tcolorbox}}
{\begin{tcolorbox}[breakable, title=Optional, colback=gray]}%
	{\end{tcolorbox}}
\usepackage{amssymb,geometry,hyperref,xcolor}
\begin{document}
	
	\maketitle
	
	\begin{abstract}
		In the paper \cite{zc} Andy Zucker and Chris Lambie-Hanson proved the consistency result
		for some coloring principle for the products of Polish spaces by at most countable many colors.
		This principle easy implies the \HLT{}. The aim of this paper is to generalize this consistency result to
		sets of colors of cardinality less than $2^{\aleph_0}$. The proof presented here differs than the proof presented in \cite{zc}.
	\end{abstract}

	\section{Required statements}
	
	\begin{te}
		Let $\lambda$ be any cardinal greater than zero.
		Then the following is true:
		\begin{enumerate}
			\item For any infinite cardinal $\mu$ holds $(2^{\mu\lambda})^\lambda=2^{\mu\lambda}>\mu$.
			\item If $\kappa$ is an infinite cardinal such that $\kappa^\lambda=\kappa$ holds, then $(\kappa^+)^\lambda=\kappa^+$ , as well as $(2^{\kappa })^\lambda=2^{\kappa}$.
		\end{enumerate}
		In particular, for any infinite cardinal $\lambda$ the class of all cardinals $\kappa$ such that $\kappa^\lambda=\kappa$ is cofinal and closed for the operations $\kappa\mapsto\kappa^+$ and $\kappa\mapsto 2^\kappa$.
	\end{te}
	\begin{dkz}
		For any infinite cardinal $\mu$
		$$
		(2^{\mu\lambda})^\lambda=2^{\mu\lambda^2}=2^{\mu\lambda}\geqslant 2^\mu>\mu
		$$
		holds.
		For any infinite cardinal $\kappa\leqslant\lambda$, it holds $\kappa^\lambda=2^\lambda>\lambda\geqslant\kappa$,
		therefore $\kappa^\lambda>\kappa$ holds. In other words, it holds
		$$
		(\forall\kappa\geqslant\aleph_0)(\kappa^\lambda=\kappa\Rightarrow\kappa>\lambda).
		$$
		Let $\kappa$ be an infinite cardinal such that $\kappa^\lambda=\kappa$;
		then $\kappa>\lambda$, so
		$$
		(2^\kappa)^\lambda=2^{\kappa\lambda}=2^\kappa
		$$
		holds.
		From $\cf(\kappa^+)=\kappa^+>\lambda$ and $\kappa^\lambda=\kappa<\kappa^+$ it follows
		$(\kappa^+)^\lambda=\kappa^+$.
	\end{dkz}
	
	For any sets $A$ and $B$ and for any infinite cardinal $\lambda$ we will denote by $\funure AB\lambda$ the set of all
	partial functions from $A$ to $B$ of cardinality less than $\lambda$.
	
	\begin{te}
		\label{deltaSistem}
		Let $\lambda$ be any infinite cardinal and let $\kappa$ and $\theta$ be regular uncountable cardinals with
		$$
		\kappa^{<\lambda}=\kappa\geqslant\lambda,\quad\theta=\kappa^+.
		$$
		Let $X$ be any set and
		let $p$ be a function that assigns to each $x\in\theta$ some partial function from $X$ to $\kappa$ of cardinality less than $\lambda$.
		Then there exists a stationary set $\theta'\subseteq\theta$, such that the set $p[\theta']$ is a $\Delta$-system.
	\end{te}
	\begin{dkz}
		Because of
		$$
		|{\textstyle\bigcup_{x\in\theta}}\dom{p(x)}|\leqslant\sum_{x\in\theta}|\dom{p(x)}|\leqslant\sum_{x\in\theta}\lambda=\theta\lambda=\theta,
		$$
		without loss of generality, we can assume that $X=\theta$.
		The set
		$$
		E'\jedef\skupsvih{\eta<\theta}{(\forall\alpha<\eta)\sup\dom{p(\alpha)}<\eta}
		$$
		is closed and unbounded as the diagonal intersection of the family
		$$
		E_\alpha\jedef\skupsvih{\eta<\theta}{\sup\dom{p(\alpha)}<\eta}
		$$
		of closed and unbounded sets. Therefore the set
		$$
		E\jedef\skupsvih{\eta\in E'}{\cf(\eta)=\kappa}
		$$
		is stationary. Let us define $f:E\longrightarrow\theta$ as $f(\alpha)=\sup(\dom{p(\alpha)}\cap\alpha)$ and
		let us choose $\alpha\in E$. It holds $\alpha\in E'$ and $\cf\alpha=\kappa$. Because of
		$|\dom{p(\alpha)}\cap\alpha|\leqslant|\dom{p(\alpha)}|<\lambda<\kappa$ it holds
		$|\dom{p(\alpha)}\cap\alpha|<\kappa$, which together with $\dom{p(\alpha)}\cap\alpha\subseteq\alpha$ and $\cf(\alpha)=\kappa$
		implies that it holds $\sup(\dom{p(\alpha)}\cap\alpha)<\alpha$, i.e. $f(\alpha)<\alpha$.
		By the pressing-down lemma, there exist
		a stationary set $S\subseteq E$ and some $\zeta<\theta$ such that
		$$
		(\forall\alpha\in S)\sup(\dom{p(\alpha)}\cap\alpha)<\zeta
		$$
		holds.
		Because of $\zeta<\theta=\kappa^+$, we can conclude that $|\zeta|\leqslant\kappa$ holds. Without loss of generality, we can assume that $\cf(\zeta)=\kappa$.
		Let $F=\funure\zeta\kappa\lambda$.
		For every set $u$ with $|u|<\lambda$ it holds
		$$
		|\kappa^u|\leqslant\kappa^{<\lambda}=\kappa.
		$$
		
		Therefore it holds
		$$
		|\skupsvih{u\subseteq\zeta}{|u|<\lambda}|=|\skupsvih{u\subseteq\kappa}{|u|<\lambda}|
		\leqslant\sum_{\mu<\lambda}\kappa^\mu\leqslant\lambda\kappa=\kappa.
		$$
		and
		$$
		|F|=|\bigcup\skupsvih{\kappa^u}{u\subseteq\zeta\land|u|<\lambda}|\leqslant\kappa^2=\kappa<\theta.
		$$
		
		For each $r\in F$ let us define a set
		$$
		A_r\jedef\skupsvih{\xi\in S}{p(\xi)\upharpoonright\xi=r}.
		$$
		
		Obviously $S=\bigcup_{r\in F}A_r$ holds, so by the regularity and uncountability of cardinal
		$\theta$, we can choose some $r\in F$ so that the set $A_r$ is stationary. Let $\theta'\jedef A_r$.
		
		Let $\xi,\eta\in\theta'$ be such that $\xi<\eta$ and let $\alpha\in\dom{p(\xi)}\cap\dom{p(\eta)} $. From
		$$
		\xi<\eta\in\theta'=A_r\subseteq S\subseteq E\subseteq E',
		$$
		it follows that $\xi<\eta\in E'$ and therefore $\alpha\leqslant\sup{\dom{p(\xi)}}<\eta$, which implies $\alpha<\eta$, which implies
		$\alpha\in\dom{p(\eta)}\cap\eta$. It follows from $\xi,\eta\in A_r$ that
		$\dom{p(\eta)}\cap\eta=\dom r=\dom{p(\xi)}\cap\xi$.
		Therefore $\alpha\in\dom r$ and $\alpha\in\dom{p_\xi}\cap\xi$ hold.
		Therefore and because of $\xi,\eta\in A_r$ it holds $p(\xi)(\alpha)=r(\alpha)=p(\eta)(\alpha)$.
		Therefore $p(\xi)\cap\ p(\eta)=r$ holds, so the set of $\skupsvih{p(\xi)}{\xi\in\theta'}$ is a $\Delta$-system with a root $r$.
	\end{dkz}
	
	\begin{te}
		\label{mnogoKonstantno}
		Let $\kappa$ be any cardinal and let $\kappa_1,\dots,\kappa_d$ be
		uncountable regular cardinals
		with $\kappa_1>\kappa$ and $\kappa_{i+1}>2^{\kappa_i}$ for all $i<d$. Then, for any $f:\kappa_1\times\dots\times\kappa_d\longrightarrow\kappa$
		there are stationary subsets $\kappa_1'\subseteq\kappa_1,\dots,\kappa_d'\subseteq\kappa_d$
		and such that the function $f$ is constant on the set $\kappa_1'\times\cdots\times\kappa_d'$.
	\end{te}
	\begin{dkz}
		We prove the statement by induction on $d$. For $d=1$ the theorem statement follows from $\cf{\kappa_1}=\kappa_1>\kappa$.
		Assume $d\geqslant 2$ and prove the inductive step. Let $f:\kappa_1\times\cdots\times\kappa_d\longrightarrow\kappa$ be arbitrary function.
		Let us define the function $G$ as follows:
		$$
		G:\kappa_d\longrightarrow \kappa^{\kappa_1\times\cdots\times\kappa_{d-1}},\quad
		(\forall a_1\in\kappa_1,\dots,a_d\in\kappa_d)G(a_d)(a_1,\dots,a_{d-1})=f(a_1,\dots,a_d).
		$$
		
		It holds
		$$
		\label{kardinalZaSledeci}
		|\kappa^{\kappa_1\times\cdots\times\kappa_{d-1}}|=\kappa^{\kappa_{d-1}}\leqslant
		(2^{\kappa})^{\kappa_{d-1}}=2^{\kappa_{d-1}}<\kappa_d,
		$$
		so there exists stationary set $\kappa_d'\subseteq\kappa_d$ such that
		$G$ is constant on the set $\kappa_d'$.
		
		Let $g$ be a function with $G(a)=g$ for all $a\in\kappa_d'$. By inductive assumption there exist
		stationary sets $\kappa_1'\subseteq\kappa_1,\dots,\kappa_{d-1}'\subseteq\kappa_{d-1}$ such that
		the function $g$ is constantly equal to some value $b$
		on the set $\kappa_1'\times\cdots\times\kappa_{d-1}'$.
		For any $a_1\in\kappa_1',\dots,a_d\in\kappa_d'$
		$$
		f(a_1,\dots,a_d)=G(a_d)(a_1,\dots,a_{d-1})=g(a_1,\dots,a_{d-1})=b
		$$
		holds.
	\end{dkz}
	
	\begin{te}
		\label{unijaFunkcija}
		Assume that $\kappa$ is regular cardinal and $\lambda$ is infinite cardinal with $\kappa^{<\lambda}=\kappa\geqslant\lambda$.
		Let $P$ denote the set of all partial functions from $\kappa$
		to $\kappa$ of cardinality less than $\lambda$. Let $\kappa_1,\dots,\kappa_d$ be arbitrary
		uncountable regular cardinals with
		\begin{enumerate}
			\item $\kappa_1=\kappa^+$,
			\item for every $i<d$ there is some regular cardinal $\mu_i\geqslant 2^{\kappa_i}$ such that $\kappa_{i+1}=\mu_i^+$ and $\mu_i^{\kappa_i}=\mu_i$ hold.
		\end{enumerate}
	
		Then for each $p:\kappa_1\times\cdots\times\kappa_d\longrightarrow P$ there are stationary sets
		$\kappa_1'\subseteq\kappa_1,\dots,\kappa_d'\subseteq\kappa_d$ such that
		$\bigcup p[\kappa_1'\times\cdots\times\kappa_d']$ is a function.
	\end{te}
	\begin{dkz}
		We prove the claim by induction on $d$. For $d=1$, the statement follows from the theorem \ref{deltaSistem}.
		Assume that $d\geqslant 2$.
		
		The set of all subsets of the set $\kappa$ of cardinality less than $\lambda$ has no more than $\kappa^{<\lambda}=\kappa$
		elements.
		There are no more mappings of any of those subsets to the set $\kappa$ than $\kappa^{<\lambda}=\kappa$.
		Therefore $|P|=\kappa$ holds.
		Let us define $F:\kappa_d\longrightarrow P^{\kappa_1\times\dots\times\kappa_{d-1}}$ as follows:
		$$
		F(x_d)(x_1,\dots,x_{d-1})=p(x_1,\dots,x_d).
		$$
		
		By inductive hypothesis, to each $x\in\kappa_d$ we can correspond stationary sets $Z_i(x)\subseteq\kappa_i$ for $i<d$
		such that the set $\bigcup F(x)[Z_1(x)\times\cdots\times Z_{d-1}(x)]$ is a function, which we will denote by $q(x)$. Moreover, for each $x_1,\dots,x_d$
		it holds $|\dom{p(x_1,\dots,x_d)}|<\lambda$, so
		$$
		\begin{array}{rcl}
			|\dom{q(x)}|&=&|\dom{{\textstyle\bigcup}(F(x)[Z_1(x)\times\cdots\times Z_{d-1}(x)])}|\leqslant\lambda|Z_1(x)|\cdots|Z_{d-1}(x)|
			\vspace{0.5em}\\
			&=&\lambda\kappa_1\cdots\kappa_{d-1}=\kappa_{d-1}.
		\end{array}
		$$

		By the theorem  \ref{deltaSistem}, there exists stationary $Z_d\subseteq\kappa_d$,
		where the set $f=\bigcup_{x\in Z_d}q(x)$ is a function.
		Therefore, for every $x_d\in Z_d$
		\begin{equation}
			\label{pSadrzanoUF}
			\arraycolsep 2pt
			\begin{array}{rcl}
				(\forall x_d\in Z_d)(\forall x_1\in Z_1(x_d),\dots,x_{d-1}\in Z_{d-1}(x_d))
				p(x_1,\dots,x_d)&=&F(x_d)(x_1,\dots,x_{d-1})
				\vspace{0.5em}\\
				&\subseteq& q(x_d)\subseteq f
			\end{array}
		\end{equation}
		holds.
		Let us define the function $G$ as follows:
		$$
		\dom G=Z_d,\quad G(x)=\torka{Z_1(x),\dots,Z_{d-1}(x)}.
		$$
		
		Because of $Z_d$ is stationary in the regular uncountable cardinal $\kappa_d$
		and $G$ maps the set $Z_d$ to the set $\Part{\kappa_1}\times\cdots\times\Part{\kappa_d}$ of the cardinality less than $\kappa_d$,
		there exists a set $\kappa_d'\subseteq Z_d$ which is stationary in $\kappa_d$ such that $G$ is a constant function on $\kappa_d'$.
		Let $\kappa_1',\dots,\kappa_{d-1}'$ be such that
		$G(x)=\torka{\kappa_1',\dots,\kappa_{d-1}'}$ holds for all $x\in\kappa_d'$.
		Let $x_d\in \kappa_d'$ be arbitrary.
		Then
		$$
		\torka{x_1,\dots,x_{d-1}}\in\kappa_1'\times\cdots\times\kappa_{d-1}'=Z_1(x_d)\times\cdots\times Z_{d-1}(x_d),
		$$
		holds for all $x_1\in\kappa_1',\dots,x_d\in\kappa_d'$, so by
		(\ref{pSadrzanoUF})
		$$
		(\forall x_1\in\kappa_1',\dots,x_d\in\kappa_d')
		p(x_1,\dots,x_d)\subseteq f
		$$
		holds.
	\end{dkz}
	
	\subsection{Trees and topology}
		
	A tree is partial ordering with the smallest element, where for every $a$ from the domain the set of all elements smaller than $a$ is well ordered.
	
	Branches are maximal chains. Here we will consider trees where each node has at least one,
	but finitely many immediate successors, and where each node has finitely many predecessors.
	Also, we limit our considerations%
\footnote{This assumption is needless, but the general case is reducible to that case.}
	to the case where each branch has
	infinitely many nodes that have more than one immediate successor. The set of all branches of the tree $T$ we denote by $[T]$.
	
	If $F$ is ultrafilter over $\omega$ and $f$ is coloring of set of nodes of the tree $T$ by finitely many colors, then we can
	define coloring $f_F$ of the set $[T]$ by the same set of colors as follows:
	$$
	f_F(b)=c\,\ekdef\,\skupsvih{n\in\omega}{f(b(n))=c}\in F.
	$$
	
	If the ultrafilter $F$ is non-principal and colorings $f$ and $g$ of the set of nodes of the tree $T$ by finitely many colors
	differs at finitely many points, then $f_F=g_F$.
	
	For any node $t$ of the tree $T$ the set of all branches $b$ of the tree $T$ such that $t\in b$ we denote by $N_t$.
	We define topology on the set $[T]$ by choosing the sets $N_t$ for all nodes $t$ of the tree $T$ as base open (or closed) sets.
	This topology is always homeomorphic to the Cantor's space. This gives us uniformization of all trees.
	
	The height of a node is the number of elements smaller than that node. The $m$-th layer in the tree $T$, denoted by $T(m)$, is the set of all its nodes of height $m$.
	
	Let $T_1,\dots,T_d$ be a sequence of trees as ordered structures. Their product is the tree denoted by $T_1\otimes\cdots\otimes T_d$ with set of nodes
	$$
	\bigcup_{m=0}^\infty(T_1(m)\times\cdots\times T_d(m)).
	$$
	and ordering defined as follows
	$$
	(a_1,\dots,a_d)<(b_1,\dots,b_d)\,\ekdef\,\bigwedge_{i=1}^da_i<b_i.
	$$
	It holds that
	$$
	(T_1\otimes\cdots\otimes T_d)(m)=T_1(m)\times\cdots\times T_d(m),\quad m\in\omega,
	$$
	and
	$$
	[T_1\otimes\cdots\otimes T_d]=[T_1]\times\cdots\times[T_d].
	$$
	
	Also, the topology of $[T_1\otimes\cdots\otimes T_d]$ is the product of the spaces $[T_1],\dots,[T_d]$.
	If $f$ is coloring of the set of the nodes of the tree $T_1\otimes\cdots\otimes T_d$ by finitely many colors and $F$ is ultrafilter over $\omega$,
	then it holds
	$$
	f_F(b_1,\dots,b_d)=c\Leftrightarrow\skupsvih{n\in\omega}{f(b_1(n),\dots,b_d(n))=c}\in F.
	$$
		
	Reading of \cite{zc} and \cite{z} is recommended because of some related results.
	
	\subsection{Polish spaces and perfect spaces}

	Polish space is complete separable metric space. Perfect space is Polish space without isolated points.
	
	Recently, the same way to generalization of results from topology of trees (Cantor's space) to arbitrary
	Polish spaces are independently used and remarked by Andy Zucker and Chris Lambie-Hanson in \cite{zc}.
	This is the following theorem and we give it with the proof.
	
	\begin{te}
		\label{gusto}
		Every perfect space has a dense $G_\delta$ subspace homeomorphic to the space ${}^\omega\omega$.
	\end{te}
	\begin{dkz}
		Let $P=\upar Xd$ is perfect space.
		Using countable base we can choose the infinite sequence $\langle U_n\,:\,n\in\omega\rangle$ of the nonempty
		open sets of diameter less than 1, where $\bigcup_{n\in\omega}U_n$ is dense and for every distinct $m,n$ it holds
		$d(U_m,U_n)>0$. Let us define the map $f$ of the set of all finite sequences of members of $\omega$ to the set of open sets of the space $P$.
		
		If $s$ is empty sequence, then $f(s)=X$. If $s$ is sequence of the length 1 and $m\in\omega$ is the only member of $s$,
		then $f(s)=U_m$.
		
		Let us assume that $s$ is some finite sequence of the members of $\omega$ such that $f(s)$ is defined and equal $V$.
		Let us denote the length of the sequence $s$ by $l$. Using the countable base we can choose some infinite
		sequence $\langle W_n\,:\,n\in\omega\rangle$ of the distinct nonempty open subsets of $V$ of diameter less than $1/(l+1)$
		such that $\bigcup_{n\in\omega}W_n$ is dense in $V$, for every distinct $m,n\in\omega$ it holds $d(W_m,W_n)>0$
		and for every $n\in\omega$ it holds $d(W_n,X\setminus V)>0$. Then for all $m\in\omega$
		we define $f(s^{\frown}m)=W_m$.
		
		Let as denote the set of all sequences of the members of $\omega$ of the length $l$ by $S_l$, the set $\bigcup f[S_l]$
		by $G_l$ and the set $\bigcap_{l\in\omega}G_l$ as $G$. The set $G$ is dense $G_\delta$ set.
		
		For any $s\in{}^\omega\omega$ there is $a\in G$ such that $\bigcap_{l\in\omega}f(s\upharpoonright l)=\{a\}$.
		This way we assign some point of $G$ to any member of ${}^\omega\omega$. This is required homeomorphism.
	\end{dkz}

	\subsection{Colorings of products of Polish spaces}
	
	Let $\kappa$ be any cardinal. By $PG(\kappa)$ we denote the following statement:
	
	\begin{center}
		\fbox{\begin{minipage}{80ex}
				\begin{center}
					\itshape
					For any natural number $d$ and Polish spaces $P_1,\dots,P_d$ and any $\alpha<\kappa$
					and $f:\dom{P_1}\times\cdots\times\dom{P_d}\longrightarrow\alpha$
					there are sets $D_1,\dots,D_d$ such that $D_i$ is somewhere dense set in the space $P_i$ for every $i$
					and where $f$ is constant on the set $D_1\times\cdots\times D_d$.
				\end{center}
		\end{minipage}}
	\end{center}

	Andy Zucker and Chris Lambie-Hanson proved in \cite{zc} that after adding at least $\beth_\omega$
	Cohen's reals to any model of ZFC, the generic extension satisfies $PG(\aleph_1)$.
	In notation used in \cite{zc} it is denoted by $PG(\aleph_0)$. The aim of this paper is to generalize this
	result for uncountable sets of colors.
	
	If $P_i$ has an isolated point $a$ then $\{a\}$ is somewhere dense set in $P_i$ and the principle is reduced to the
	product of other spaces. Therefore, the general case is reduced to products of perfect spaces.
	According to the theorem \ref{gusto} the general case is reduced to the case that all of spaces
	$P_1,\dots,P_d$ are ${}^\omega\omega$.
	
	Andy Zucker and Chris Lambie-Hanson derived Halpern and L\" auchli theorem from this consistency fact in \cite{zc}.
	For completeness of this paper, we give the derivation here.
	
	For the first, we will derive Halpern and L\" auchli theorem from $PG(\aleph_0)$. Let
	$T_1,\dots,T_d$ be any trees and $f$ be any coloring $f$ the set of nodes of the tree $T_1\otimes\cdots\otimes T_d$ by
	the finitely many colors. Let us choose any non-principal ultrafilter $F$ over $\omega$.
	Then $f_F$ is coloring of the set $[T_1]\times\cdots\times[T_d]$ by the same set of colors.
	
	According to $PG(\aleph_0)$, there are somewhere dense sets $S_1,\dots,S_d$ of the spaces $[T_1],\dots,[T_d]$
	such that $f_F$ is equal to some constant $c$ on the set $S_1\times\cdots\times S_d$.
	Let us choose the nodes $t_1,\dots,t_d$ such that $S_i$ is dense in $N_{t_i}$ in the space $[T_i]$
	for all $i$. We can choose nodes $t_1,\dots,t_d$ of the trees $T_1,\dots,T_d$ at the same height $l$.
	
	Let us denote all distinct children of the node $t_i$ in the tree $T_i$ by $s^i_1,\dots,s^i_{n_i}$
	and let us choose $b^i_j\in S_i$ such that $s^i_j\in b^i_j$. For every choice of
	$j_1,\dots,j_d\in\omega$ such that $\bigwedge_{i=1}^d0<j_i\leqslant n_i$ it holds
	$f_F(b^1_{j_1},\dots,b^d_{j_d})=c$ and therefore the set
	\begin{equation}
		\label{uFilteru}
		\skupsvih{n\in\omega}{f(b^1_{j_1}(n),\dots,b^d_{j_d}(n))=c}
	\end{equation}
	is in $F$. Therefore, intersection $I$ of sets of form (\ref{uFilteru})  for all finitely many choices
	of $j_1,\dots,j_d$ is in $F$. Because of $F$ is non-principal ultrafilter, we can choose some $n\in I$ such that
	$n>l+1$ holds. Let us define $D_i$ as the set $\{b^i_1(n),\dots,b^i_d(n)\}$.
	
	Every child of $t_i$ in the tree $T_i$ has successor in the set $D_i$ and $f$ is equal to the constant $c$ on the
	set $D_1\times\cdots\times D_d$.
	
	Let us prove the Halpern and L\" auchli theorem in $\ZFC$. Choose we some trees $T_1,\dots,T_d$ and coloring $f$.
 	These object have isomorphic image such that it's transitive closure is countable.
	Let us choose some countable transitive model $M$ for adequate finite fragment of $\ZFC$ containing these objects and
	make the forcing extension $M[G]$ where $PG(\aleph_0)$ holds.
	This model contains witnesses $t_1,\dots,t_d$ and $D_1,\dots,D_d$ in $M[G]$
	and validation is absolute for transitive models. Therefore, witnesses are also valid in the universe.

	\subsection{Consistency of PG with ZFC}
	
	Let us prove that the $PG$ principle also holds in an important model that satisfies the $\ZFC$ axioms.
	Before that, we will introduce the notion of $\Delta^d$-system and prove the generalization of the theorem
	\ref{unijaFunkcija}.
	
	Thereafter $x_1,\dots,x_d$ will be denoted by $\bar x$ and $\torka{x_1,\dots,x_d}$ will be denoted by $\vec x$.
	
	\begin{de}
		Let $X_1,\dots,X_d,D$ and $Y$ be the sets, $p$ the mapping of the set $X_1\times\cdots\times X_d$ into the set of partial functions from $D$ to $Y$ and $r_1,\dots,r_d$ partial functions from $D$ to $Y$ such that the following holds:
		\begin{enumerate}
			\item $r_1\subseteq\cdots\subseteq r_d$,
			\item for each $k<d$ and all $x_1\in X_1,\dots,x_d\in X_d$ the set
			$$
			\skupsvih{p(x_1,\dots,x_d)\cap r_{k+1}}{x_1\in X_1,\dots,x_d\in X_d}
			$$
			forms a $\Delta$-system with the root $r_k$,
			\item for each $k<d$ the value of $p(x_1,\dots,x_d)\cap r_{k+1}$ depends only on $x_1,\dots,x_k$,
			\item for each $q$ there are sets $F_1,\dots,F_d$ of cardinality not greater than $|q|$ such that the following holds:
			\begin{enumerate}
				\item
				\label{prvok}
				The statement $\dom{p(\bar x)}\cap q\subseteq\dom{r_1}$ holds for every $\vec x\in(X_1'\setminus F_1)\times\cdots\times (X_d'\setminus F_d)$.
				\item
				\label{drugok}
				The statement $\dom{p(\bar x)}\cap q\subseteq\dom{r_{k+1}}$ holds
				for each $k<d$ and any
				$$
				\vec x\in X_1'\times\cdots\times X_k'\times(X_{k+1}'\setminus F_{k+1})\times\cdots\times (X_d'\setminus F_d).
				$$
			\end{enumerate}
		\end{enumerate}
		then we say that $\torka{p,X_1,\dots,X_d,r_1,\dots,r_d}$ is a $\Delta^d$-system.
	\end{de}
	
	Obviously, the notion of a $\Delta^1$-system coincides with the notion of a $\Delta$-system.
	
	\begin{te}
		\label{proizvodDelta}
		Let $\kappa$ be regular infinite cardinal and $\lambda$ be any infinite cardinal where $\kappa^{<\lambda}=\kappa\geqslant\lambda$ holds
		and let $\kappa_1,\dots,\kappa_d$ be cardinals such that that $\kappa_1=\kappa^+$ and that for every $i<d$ there is some $\mu_i\geqslant 2^{\kappa_i}$ such that 
		$\kappa_{i+1}=\mu_i^+$ and $\mu_i^{\kappa_i}=\mu_i$ hold. Let $X_1,\dots,X_d$ and $D$ are sets such that
		$\bigwedge_{i=1}^d|X_i|=\kappa_i$, $D\supseteq\bigcup_{i=1}^d$ and $|D|=\kappa_d$ hold.
		Let $p:X_1\times\cdots\times X_d\longrightarrow\funure D\kappa\lambda$ be such that
		$$
		(\forall x_1\in X_1,\dots,x_d\in X_d)\{x_1,\dots,x_d\}\subseteq\dom{p(x_1,\dots,x_d)}
		$$
		hold.
		Then there exist sets $X_1',\dots,X_d'$ and $r_1,\dots,r_d$ such that the following holds:
		\begin{enumerate}
			\item The set $X_i'$ is the subset of $X_i$ of the cardinality%
\footnote{Formulation with stationary sets is also valid, but needless for this paper.}
			 $\kappa_i$ for every $i$.
			\item $\torka{p,X_1',\dots,X_d',r_1,\dots,r_d}$ is a $\Delta^d$-system,
			\item $|r_1|<\lambda$ and $\bigwedge_{i=2}^d|r_i|\leqslant\kappa_{d-1}$.
		\end{enumerate}
	\end{te}
	\begin{dkz}
		By the theorem \ref{unijaFunkcija}, without loss of generality we can assume that $\bigcup p[X_1\times\cdots\times X_d]$ is a function.
		
		Let $d=1$. By the theorem \ref{deltaSistem}, there exists $X_1'\subseteq X_1$ such that $|X_1'|=\kappa_1$ holds and
		the set $p[X_1']$ forms a $\Delta$-system with root $r_1$.
		Obviously $|r_1|<\lambda$ and $|\bigcup p[X_1']|\leqslant\kappa_1$.
		
		Let $d\geqslant 2$.
		For each $x_d\in X_d$ let us define the function
		$$
		p'(x_d)=\bigcup\skupsvih{p(x_1,\dots,x_d)}{x_1\in X_1,\dots,x_{d-1}\in X_{d-1}}.
		$$
		
		By the theorem \ref{deltaSistem}, there exists $X\subseteq X_d$ such that $|X|=\kappa_d$
		and that the set $p'[X]$ forms a $\Delta$-system with some root $r_d$.
		It holds
		\begin{equation}
			\label{kardinalnostPPrim}
			(\forall x\in X)|p'(x)|\leqslant\kappa_{d-1}
		\end{equation}
		which implies $|r_d|\leqslant\kappa_{d-1}$. For every $x_1\in X_1,\dots,x_d\in X_d$ it holds
		$$
		x_d\in X\Rightarrow\{x_1,\dots,x_d\}\subseteq\dom{p(x_1,\dots,x_d)}\subseteq\dom{p'(x_d)}.
		$$
		Therefore $x_d\in\dom{p'(x_d)}$ holds, which implies
		$$
		\bigcup\dom{p'[X]}\supseteq X.
		$$
		
		Therefore and from (\ref{kardinalnostPPrim}) it follows that $|p'[X]|=\kappa_d$ holds, therefore there exists $X'\subseteq X$
		such that $|X'|=\kappa_d$ and such that the function $p'$ is an injection on the set $X'$.
		Let us define the function $p''$ as follows:
		$$
		\dom{p''}=X',\quad
		p''(x_d)(x_1,\dots,x_{d-1})=p(x_1,\dots,x_d)\cap r_d.
		$$
		
		By the inductive hypothesis we can define sets
		$$
		X_1'(x_d),\dots,X_{d-1}'(x_d)\mbox{ and } r_1(x_d),\dots,r_{d-1}(x_d),
		$$
		for every $x_d\in X'$, so that the statement of the theorem applied to the function $p''(x_d)$ instead of the function $p$ holds.
		By the theorem \ref{mnogoKonstantno}, there exists a set $X_d'\subseteq X'$ such that $|X_d'|=\kappa_d$
		and that the functions $X_1',\dots,X_{d-1}'$, $p''$ and $r_1,\dots,r_{d-1}$ are constant on the set $X_d'$.
		The values of these functions on the set $S_d$ we will denote by the same labels used for those functions.
				
		Let $x_1\in X_1',\dots,x_d\in X_d'$. Then
		$$
		p(\bar x)\cap r_d=p''(x_d)(x_1,\dots,x_{d-1})=p''(x_d')(x_1,\dots,x_{d-1})=p(x_1,\dots,x_{d-1},x_d)\cap r_d
		$$
		holds. Therefore the following formula
		$$
		(\forall x_1\in X_1',\dots,x_d\in X_d')\Big(\bigwedge_{i=1}^kx_i=y_i\Rightarrow p(\bar x)\cap r_{k+1}=p(\bar y)\cap r_{k+1}\Big)
		$$
		holds for $k=d-1$. Let us show that it also holds for $k<d-1$.

		Let $k<d-1$ and let $x_1,y_1\in X_1',\dots, x_d,y_d\in X_d'$ be such that $\bigwedge_{i=1}^kx_i=y_i$ holds.
		Then the following
		$$
		\begin{array}{rcl}
			p(\bar x)\cap r_{k+1}&=&p(\bar x)\cap r_d\cap r_{k+1}=p''(x_d)(x_1,\dots,x_{d-1})\cap r_{k+1}
			\vspace{0.5em}\\
			&=&
			p''(y_d)(x_1,\dots,x_{d-1})\cap r_{k+1}=
			p''(y_d)(y_1,\dots,y_{d-1})\cap r_{k+1}
			\vspace{0.5em}\\
			&=&
			p(\bar y)\cap r_d\cap r_{k+1}=p(\bar y)\cap r_{k+1}
		\end{array}
		$$
		holds.
		Let us choose arbitrary $\vec x\in X_1'\times\cdots\times X_d'$.
		Since the set $p'[X_d']$ forms a $\Delta$-system with the root $r_d$,
		it is true that $p(\bar x)$ and $r_d$ are restrictions of the same function $p'(x_d)$. Therefore
		$$
		\dom{p(\bar x)\cap r_d}=\dom{p(\bar x)}\cap\dom{r_d}
		$$
		holds.
		Also, $r_{k+1}\subseteq r_d$ holds for every $k<d$, so $p(\bar x)$ and $r_{k+1}$ are also restrictions of the same function $p'(x_d )$, so
		$$
		\dom{p(\bar x)\cap r_{k+1}}=\dom{p(\bar x)}\cap\dom{r_{k+1}}
		$$
		holds.
		Let $q\subseteq D$ be arbitrary. Since the function $p'$ is
		injection on the set $X_d'$ and the set $p'[X_d']$ forms a $\Delta$-system with root $r_d$, there exists $F_d\subseteq X_d'$ such that $|F_d|\leqslant|q| $ and
		\begin{equation}
			\label{rd}
			(\forall\vec x\in X_1'\times\cdots\times X_{d-1}'\times(X_d'\setminus F_d))\dom{p(\bar x)}\cap q\subseteq\dom{r_d}
		\end{equation}
		hold.
		By the inductive hypothesis,
		there are subsets $F_1,\dots,F_{d-1}$ of the sets $X_1',\dots,X_{d-1}'$, having cardinalities no greater than $|q|$,
		such that
		\begin{equation}
			\label{r1}
			(\forall\vec x\in(X_1'\setminus F_1)\times\cdots\times(X_d'\setminus F_d))\dom{p(\bar x)}\cap\dom{r_d}\cap q\subseteq\dom{r_1},
		\end{equation}
		holds, as well as that
		\begin{equation}
			\label{rk1}
			(\forall\vec x\in X_1'\times\cdots\times X_k'\times(X_{k+1}'\setminus F_{k+1})\times\cdots\times(X_d'\setminus F_d))
			\dom{p(\bar x)}\cap\dom{r_d}\cap q\subseteq\dom{r_{k+1}}
		\end{equation}
		holds for every $k<d-1$.
		The formula (\ref{rk1}) obviously holds for $k=d-1$.
		Assume that
		$$
		\vec x\in(X_1'\setminus F_1)\times\cdots\times(X_d'\setminus F_d)
		$$
		holds. Then, by (\ref{rd}) and (\ref{r1}) $\dom{p(\bar x)}\cap q\subseteq\dom{r_1}$ holds.
		Now suppose that
		$$
		k<d\mbox{ and }\vec x\in X_1'\times\cdots\times X_k'\times(X_{k+1}'\setminus F_{k+1})\times\cdots\times(X_d'\setminus F_d)
		$$
		holds.
		Then by (\ref{rd}) and (\ref{rk1}),  $\dom{p(\bar x)}\cap q\subseteq\dom{r_{k+1}}$ holds.
	\end{dkz}
	\vspace{0.5em}
	
	The previous theorem also has a finitary variant, which we do not give here because it is irrelevant for the set theory.
	Note that forcing with finite partial functions from $\omega$ to $\omega$ ordered by reverse inclusion adds one Cohen's real as
	forcing by countable atomless poset.
	
	\begin{te} \label{malaMG} Let $M$ be a countable transitive model. Let $\theta_0$ be uncountable regular cardinal in $M$
		and for every $i\in\omega$ let $\theta_{i+1}$ be $2^{\theta_i}$ in $M$. Let $\lambda=\sup_i\theta_i$ and
		$M[G]$ is generic extension obtained by adding at least $\lambda$
		Cohen's reals. Then $M[G]\models PG(\lambda)$ holds.
	\end{te}
	\begin{dkz}
		Let us define poset $P$ as finite partial functions from $\lambda$ to set of all finite sequences of elements of $\omega$
		with the following ordering:
		$$
		p\leqslant q\,\ekdef\,\dom q\subseteq\dom p\,\land\,(\forall x\in\dom q)q(x)\subseteq p(x).
		$$
		
		A generic object will represent a family of members of ${}^\omega\omega$ indexed by ordinals from the set $\lambda$.
		Let us define the name
		$$
		\dot b\jedef\skupsvih{\upar q{\upar\sigma{\upar nk}\check{}\,}}{\sigma\in\dom q\land n\in\dom{q(\sigma)}\land q(\sigma)(n)=k}.
		$$
		
		In other words, $\dot b$ is the name for the function that maps every $\sigma\in\lambda$ to generic sequence with index $\sigma$.
		Let's assume that $M[G]\models\lnot PG(\lambda)$. Then there are $p_0\in G$, name $\dot f$, $d\in\omega\setminus\{0\}$
		and $\mu<\lambda$ such that it holds that
		$$
		M\models\,p_0\Vdash\dot f:({}^\omega\omega)^d\longrightarrow\mu
		$$
		and in $M$ the condition $p_0$ forces that there are no somewhere dense sets $D_1,\dots,D_d$ in the space ${}^\omega\omega$
		such that $\dot f$ is constant on $D_1\times\cdots\times D_d$. Let us choose $j\in\omega$ such that $\theta_j>\mu$ holds.
		Let $\kappa_0$ be $\theta_i$
		and $\kappa_{i+1}$ be $(2^{\kappa_i})^{++}$ for every $i\in\omega$. For every cardinal $\kappa$ it holds
		$\kappa<2^\kappa<(2^\kappa)^{++}\leqslant 2^{2^{2^\kappa}}$. Therefore $\sup_i\kappa_i=\lambda$ holds.

		To each $\vec\sigma\in\kappa_1\times\cdots\times\kappa_d$ we associate the condition
		$p(\bar\sigma)\leqslant p_1$ and some $k(\bar\sigma)\in\mu$ such that
		$$
		M\models p(\bar\sigma)\Vdash_P\dot f(\dot b(\sigma_1),\dots,\dot b(\sigma_d))
		=k(\bar\sigma)\,,\quad\{\sigma_1,\dots,\sigma_d\}\subseteq\dom{p(\bar\sigma)},
		$$
		holds, where the condition $p(\bar\sigma)$
		maps all the elements of its domain into strings of the same length
		$l(\bar\sigma)\in\omega$.
		
		For each $i\leqslant d$, we define a function $c_i$ whose domain is $\kappa_1\times\cdots\times\kappa_d$,
		such that for every $\vec\sigma$ from the domain
		$$
		c_i(\bar\sigma)=p(\bar\sigma)(\sigma_i)
		$$
		holds.
		In other words, $c_i(\bar\sigma)$ is the information carried by the condition $p(\bar\sigma)$
		about the generic branch at position $\sigma_i$. Let us define a function $g$ with the same domain such that
		$$
		g(\bar\sigma)=\torka{k(\bar\sigma),l(\bar\sigma),c_1(\bar\sigma),\dots,c_d(\bar\sigma)}
		$$
		holds.
		Note that $|\mathrm{range}(g)|^M<\kappa_0$. By the theorems
		\ref{mnogoKonstantno} and \ref{unijaFunkcija}, there are in $M$ infinite subsets
		$H_1,\dots,H_d$
		of sets $\kappa_1,\dots,\kappa_d$ such that for all $i$
		it holds $|H_i|^M=\kappa_i$, that the function $g$
		on the set $H_1\times\cdots\times H_d$ is constantly equal to some $\torka{k,l,c_1,\dots,c_d}$ and that
		all conditions from the set $p[H_1\times\cdots\times H_d]$ are compatible. We will denote $H_1\times\cdots\times H_d$
		by $H$.
		
		By the theorem \ref{proizvodDelta}, the sets $H_1,\dots,H_d$
		could be chosen so that there are $r_1,\dots,r_d$ from $M$ such that the following holds:
		\begin{enumerate}
			\item $\torka{p,H_1,\dots,H_d,r_1,\dots,r_d}$ is a $\Delta^d$-system,
			\item $|r_1|<\aleph_0$ and $\bigwedge_{i=2}^d|r_i|^M=\kappa_{i-1}$.
		\end{enumerate}
		
		With the symbolic from the definition of the $\Delta^d$-system, we will denote by $F_1(q_0),\dots,F_d(q_0)$ the sets $F_1,\dots,F_d$ corresponding to the set $q=\dom{q_0}$
		for a given condition $q_0$. Let's define
		$$
		\dot T\jedef\skupsvih{\upar{p(\bar\sigma)}{\torka{\sigma_1,\dots,\sigma_d}\,{}\check{}\,}}{
			\vec\sigma\in H}.
		$$
		
		Let us prove the following:
		$$
		M\models r_1\Vdash_P(\forall\vec\sigma\in\dot T)\dot f(\dot b(\sigma_1),\dots,\dot b(\sigma_d))=k.
		$$
		
		Otherwise, there exist $\vec\sigma\in H$ and $q_0\leqslant r_1$ such that
		$$
		M\models q_0\Vdash_P\dot f(\dot b(\sigma_1),\dots,\dot b(\sigma_d))\neq k,
		$$
		$$
		M\models q_0\Vdash_P\vec\sigma\in\dot T
		$$
		holds. The last formula means that $q_0\leqslant p(\bar\sigma)$ and therefore
		$$
		M\models q_0\Vdash_P\dot f(\dot b(\sigma_1),\dots,\dot b(\sigma_d))=k,
		$$
		which is a contradiction. Let us define
		$$
		[s]\jedef\skupsvih{x\in {}^\omega\omega}{s\subseteq x},\quad s\mbox{ is a finite sequence of members of }\omega.
		$$
		
		Let us choose any infinite sets $I_1,\dots,I_d\in M$ such that $I_k\subseteq H_i$ holds for every $i$.
		Let us denote $I_1\times\cdots\times I_d$ by $I$.
		Let $s_1,\dots,s_d$ are some finite sequences of the members of $\omega$ such that $s_1\supseteq c_1,\dots,s_d\supseteq c_d$ holds.
		Let us prove the following:
		\begin{equation}
			\label{rForsira1}
			M\models r_1\Vdash_P(\exists\vec\sigma\in I\cap\dot T)\torka{\dot b(\sigma_1),\dots,\dot b(\sigma_d)}\in[s_1]\times\cdots\times[s_d].
		\end{equation}
		
		Otherwise, there exist $q_0\leqslant r_1$ such that
		$$
		M\models q_0\Vdash_P\lnot(\exists\vec\sigma\in I\cap\dot T)\torka{\dot b(\sigma_1),\dots,\dot b(\sigma_d)}\in[s_1]\times\cdots\times[s_d].
		$$
		holds. Let us choose any $\vec\tau$ from the set $I$ such that $\bigwedge_{i=1}^d\tau_i\not\in(F_i(q_0)\cup\dom{q_0})$ holds.
		From $\dom{p(\bar\tau)}\cap\dom{q_0}\subseteq\dom{r_1}$ and $p(\bar\tau),q_0\leqslant r_1$ we can conclude that $p(\bar\tau)||q_0$ holds.
		
		The information that the condition $p(\bar\tau)$ contains about the Cohen's
		real in place $\tau_i$ is $c_i$ for all $i$. By the choice of the element $\tau_i$, the condition $q_0$ contains no information about this Cohen's real.
		Let us denote by $q_1$ the greatest condition bellow the condition $p(\bar\tau)$ so that $q_1$ contains the information $s_i$
		about Cohen's real
		at place $\tau_i$ for all $i$. Then $q_1\leqslant p(\bar\tau)$ and $q_1||q_0$ is valid.
		Let us choose the condition $q_2$ so that $q_2\leqslant q_0,q_1$ holds. Due to $q_2\leqslant p(\bar\tau)$,
		$$
		M\models q_2\Vdash_P\torka{\tau_1,\dots,\tau_d}\in\dot T
		$$
		holds. Because of $q_2\leqslant q_1$ it holds
		$$
		M\models q_2\Vdash_P\bigwedge_{i=1}^d\dot b(\tau_i)\in[s_i],
		$$
		which contradicts $q_2\leqslant q_0$ and the choice of $q_0$. This contradiction proves
		(\ref{rForsira1}). Assume we that $d>1$. Let $s_d$ is any finite sequence of members of $\omega$ such that
		$s_d\supseteq c_d$ holds.
		Let us prove that in the model $M$ the condition $r_1$ forces that for every finite subset $F$ of the set
		$$
		\skupsvih{\torka{\sigma_1,\dots,\sigma_{d-1})}}{(\exists\sigma_d)\vec\sigma\in\dot T\cap I}
		$$
		there is some $\sigma_d\in I_d$ such that $\dot b(\sigma_d)\in[s_d]$ and $\vec\sigma\in\dot T\cap I$
		for every $\torka{\sigma_1,\dots,\sigma_{d-1}}\in F$.
		Otherwise, there is $q_0\leqslant r_1$ which forces negation of this statement with decided set $F$.
		
		Let $F$ be the set
		$$
		\{\torka{\sigma^1_1,\dots,\sigma^1_{d-1}},\dots,\torka{\sigma^m_1,\dots\sigma^1_{d-1}}\}.
		$$
		
		Without loosing of generality we can assume that for some decided $\sigma^1_d,\dots,\sigma^m_d$
		the condition $q_0$ forces that $\torka{\sigma^i_1,\dots,\sigma^i_d}\in\dot T\cap I$ for all $i$.
		Because of
		$$
		M\models q_0\Vdash_P\bigwedge_{i=1}^m\torka{\sigma^i_1,\dots,\sigma^i_d}\in I,
		$$
		it holds
		$$
		\torka{\sigma^1_1,\dots,\sigma^1_d},\dots,\torka{\sigma^m_1,\dots,\sigma^m_d}\in I.
		$$
		Because of
		$$
		M\models q_0\Vdash_P\bigwedge_{i=1}^m\torka{\sigma^i_1,\dots,\sigma^i_d}\in\dot T,
		$$
		it holds
		$$
		q_0\leqslant p(\sigma^1_1,\dots,\sigma^1_d),\dots,p(\sigma^m_1,\dots,\sigma^m_d).
		$$
		Let us choose $\sigma_d\in I_d\setminus(F_d(q_0)\cup\dom{q_0})$.
		For each $i$
		$$
		\dom{p(\sigma^i_1,\dots,\sigma^i_{d-1},\sigma_d)}\cap\dom{q_0}\subseteq\dom{r_d}
		$$
		holds and therefore
		$$
		\begin{array}{rcl}
			\dom{p(\sigma^i_1,\dots,\sigma^i_{d-1},\sigma_d)}\cap\dom{q_0}&=&\dom{p(\sigma^i_1,\dots,\sigma^i_{d-1},\sigma_d)}\cap\dom{r_d}\cap\dom{q_0}
			\vspace{0.5em}\\
			&=&\dom{p(\sigma^i_1,\dots,\sigma^i_{d-1},\sigma^i_d)}\cap\dom{r_d}\cap\dom{q_0}
		\end{array}
		$$
		which together with $q_0\leqslant p(\sigma^i_1,\dots,\sigma^i_d)$ and
		$$
		p(\sigma^i_1,\dots,\sigma^i_{d-1},\sigma_d)\upharpoonright\dom{r_d}=p(\sigma^i_1,\dots,\sigma^i_d)\upharpoonright\dom{r_d}
		$$
		implies that $q_0||p(\sigma^i_1,\dots,\sigma^i_{d-1},\sigma_d)$ holds.
		For condition $q_1$ defined as
		$$
		q_1\jedef\inf\skupsvih{p(\sigma^i_1,\dots,\sigma^i_{d-1},\sigma_d)}{i\in\{1,\dots,m\}}
		$$
		it holds $q_0||q_1$. The information that $q_1$ contains about infinite sequence at position $\sigma_d$ is $c_d$
		while from $\sigma_d\not\in\dom{q_0}$ we conclude that $q_0$ does not contain information about it.
		Therefore, there is condition $q_2\leqslant q_0,q_1$ such that information that $q_2$ contains about infinite sequence at position $\sigma_d$
		is $s_d$. Therefore,
		$$
		M\models q_2\Vdash_P(\sigma_d\in I_d\land\dot b(\sigma_d)\in[s_d]),
		$$
		and
		$$
		M\models q_2\Vdash_P\torka{\sigma^1_1,\dots,\sigma^1_{d-1},\sigma_d},\dots,\torka{\sigma^1_1,\dots,\sigma^1_{d-1},\sigma_d}\in\dot T\cap I,
		$$
		which is contradiction. Obviously, the similar statement holds if we replace sets $\kappa_1,\dots,\kappa_d$ by
		$\kappa_{n_1},\dots,\kappa_{n_d}$ for some $n_1,\dots,n_d\in\omega$ where $n_1<\cdots<n_d$ holds.
		
		Moreover, the similar statement holds if we replace the sets $\kappa_1,\dots,\kappa_d$ by sets
		$\kappa_{n_{\pi(1)}},\dots,\kappa_{n_{\sigma_d}}$ for some distinct $n_1,\dots,n_d\in\omega$ and
		permutation $\pi$ of the set $\{1,\dots,d\}$ where $n_{\sigma(1)}<\cdots<n_{\sigma(d)}$ holds. 
		
		Let us choose permutations $\pi_1,\dots,\pi_d$ of the set $\{1,\dots,d\}$ such that $\pi_i(d)=i$.
		Let us choose $n^i_j\in\omega$ for $i,j\in\{1,\dots,d\}$ such that
		$$
		\pi_i(n^i_1)<\cdots<\pi_i(n^i_d)
		$$
		holds for all $i\in\{1,\dots,d\}$ and where
		$$
		\pi_i(n^i_j)\leqslant\pi_{i+1}(n^{i+1}_j)
		$$
		holds for all $i,j\in\{1,\dots,d\}$ with $i\neq d$.
		
		We can choose in $M$ the sets $H^d_1\subseteq\kappa_{n^d_1},\dots,H^d_d\subseteq\kappa_{n^d_d}$ such that
		$\bigcup p[H^d_1\times\cdots\times H^d_d]$ is a function and where $|H^d_j|^M=\kappa_{n^d_j}$ and where we obtain some $\Delta^d$-system
		with $t_d$ as the smallest root in system. Then we obtain $c_1,\dots,c_d,l$ and $k$ and these values will not be changed later
		because of we will just shrink system.
		
		For every $i<d$ we can choose in $M$ some sets $H^i_1\subseteq H^{i+1}_1,\dots,H^i_d\subseteq H^{i+1}_d$ such that
		$|H^i_j|^M=\kappa_{n^i_j}$ for all $j$ and where we obtain some $\Delta^d$-system with $t_i$ as the smallest root in system.
		
		Of course, the roots $t_1,\dots,t_d$ are restrictions of the function $\bigcup p[H^d_1\times\cdots\times H^d_d]$
		and therefore $t\jedef t_1\cup\cdots\cup t_d$ is the function and moreover condition.
		
		Let us choose in $M$ infinite sets $I_1\subseteq H^1_1,\dots,I_d\subseteq H^1_d$ and define $I$ as
		$I_1\times\cdots\times I_d$. For each $i\in\{1,\dots,d\}$ let $\langle s^i_j\,:\,j\in\omega\rangle$ is the sequence of all finite sequences $u$ of the members of $\omega$
		such that $u\supseteq c_i$ holds.
		Let $H$ be some generic filter such that $t\in H$ holds.
		
		The following construction we perform in $M[H]$.
		For the first, we can choose $\sigma^1_0\in I_1,\dots,\sigma^d_0\in I_d$ such that
		$$
		\torka{\sigma^1_0,\dots,\sigma^d_0}\in\dot T_H
		$$
		and
		$$
		\dot b_H(\sigma^1_0)\in[s^1_0],\dots,\dot b_H(\sigma^d_0)\in[s^d_0].
		$$
		hold. Let us presume that $\sigma^1_i,\dots,\sigma^d_i$ are defined for some $i\in\omega$. We can choose
		$\sigma^1_{i+1}\in I_1$ such that $\dot b_H(\sigma^1_{i+1})\in[s^1_{i+1}]$ holds and (if $d\geqslant 2$) for all $m_2,\dots,m_d\leqslant i$ it holds
		$$
		\torka{\sigma^1_{i+1},\sigma^2_{m_2},\dots,\sigma^d_{m_d}}\in\dot T_H.
		$$
		
		Let us presume that $\sigma^1_{i+1},\dots,\sigma^j_{i+1}$ are defined for some $j<d$. We can choose
		$\sigma^{j+1}_{i+1}\in I_{j+1}$ such that $\dot b_H(\sigma^{j+1}_{i+1})\in[s^{j+1}_{i+1}]$ and (if $j+1<d$)
		for all $m_1,\dots,m_d\in\omega$ the following
		$$
		(\forall k<j+1)m_k\leqslant i+1\,\land\,m_{j+1}=i+1\,\land\,(\forall k>j+1)m_k\leqslant i
		$$
		implies that $\torka{\sigma^1_{m_1},\dots,\sigma^d_{m_d}}\in\dot T_H$ holds.
		The sets $D_1,\dots,D_d$ defined as
		$$
		D_i\jedef\skupsvih{\dot b_H(\sigma^i_j)}{j\in\omega}
		$$
		are somewhere dense and the function $\dot f_H$ is constant on the set $D_1\times\cdots\times D_d$.
		This contradicts to $p_0\in H$.
	\end{dkz}

	Therefore in the model with $\beth_{\omega_1}$ added Cohen reals it holds $PG(2^{\aleph_0})$.
	In other words the $PG$ principle is valid for less than $2^{\aleph_0}$ colors.
	Applications of the case with uncountable many colors are expected in inductive proofs and recursive definitions over the minimal well-orderings of sets of cardinality $2^{\aleph_0}$.

\end{document}